 \documentclass{article}
 \usepackage{amstex}

\newcommand\Frac{\mathop{\mathrm{Frac}}}
\newtheorem{theorem}{Theorem}
 
\newtheorem{claim}[theorem]{Main Claim}
\newenvironment{pf}{\paragraph{Proof}}{\par\medskip}             
\newenvironment{pfof}[1]{\paragraph{Proof of #1}}{\par\medskip}  

\newcommand\hA{\widehat A}

\newcommand\Ga{\Gamma}

\newcommand\wave{\widetilde}
\newcommand\aff{\Bbb A}
\newcommand\C{\Bbb C}
\newcommand\Z{\Bbb Z}
\newcommand\Spec{\operatorname{Spec}}
\newcommand\QED{\ifhmode\unskip\nobreak\fi\quad {\rm Q.E.D.}} 
\newcommand\qed{\ifhmode\unskip\nobreak\fi\quad $\Box$}       
\DeclareMathSymbol{\onto}  {\mathrel}{AMSa}{"10}

\title{Akizuki's counterexample}

\author{Miles Reid
\thanks{This paper was written during a stay at the semestre ``Surfaces de
Riemann et fibr\'es vectoriels'' of the Centre Emile Borel, paid for by the
EEC HCM project AGE (Algebraic Geometry in Europe), contract number ERBCHRXCT
940557.}}
\date{Mar 1995}

\begin{document}

\maketitle
 \begin{abstract}
Following \cite{Akizuki}, I construct a Noetherian local
integral domain $C_M$ whose normalisation (integral closure) is not finite
over $C_M$. My proof follows closely Akizuki's ingenious calculations.
 \end{abstract}

Let $A$ be a DVR with local parameter $t$ and residue field $k=A/(t)$, and
$\hA$ the completion of $A$. There are no restrictions on the characteristic of
$A$ or $k$, but I assume that $\hA$ contains a transcendental element over $A$.
(For DVRs of interest, $\hA$ usually has infinite transcendence degree over
$A$.) The rings $B,C$ constructed below, and their localisations, are
intermediate subrings between $A$ and $\hA$.

The construction depends on a power series
 \begin{equation}
z=z_0=a_0+a_1t^{n_1}+a_2t^{n_2}+\cdots\in\hA,
 \label{z=z(t)}\end{equation}
(not just on the element $z$). Assume:
 \begin{enumerate}
\renewcommand\labelenumi{(\arabic{enumi})}\setcounter{enumi}{\value{equation}}
 \item Each $a_i\in A$ is a unit.
 \item\label{ineq:nr} $n_r\ge2n_{r-1}+2$ for every $r\ge1$, where I
set $n_0=0$; for example, the smallest possible choice is
$n_r=2(2^r-1)=0,2,6,14,30,\dots$.
 \item $z$ is transcendental over $A$, so that $A\subset A[z]\subset\hA$ is a
polynomial extension.
 \setcounter{equation}{\value{enumi}}
 \end{enumerate}

Akizuki's construction is as follows: for $r\ge0$, let
 \begin{equation}
z_r={z_0-\hbox{first $r$ terms}\over
t^{n_r}}=a_r+a_{r+1}t^{m_{r+1}}+\cdots,
 \nonumber\end{equation}
where
 \begin{equation}
m_r=n_r-n_{r-1}\quad\text{so that (\ref{ineq:nr}) gives}\quad
2m_r\ge n_r+2.
 \label{m_r=}\end{equation}
Then the $z_r$ satisfy the identities
 \begin{eqnarray}
z_r-a_r&=&t^{m_{r+1}}z_{r+1},
 \label{eq:z_r=}\\
t^{n_r}z_r&=&z_0-\sum_{i=0}^{r-1} a_it^{n_i},\quad
\text{with $\sum_{i=0}^{r-1} a_it^{n_i}\in A$.}
 \label{eq:tnrz_r=}\end{eqnarray}

Then set $B=A[z_0,z_1,\dots]=A[(z_0-a_0),(z_1-a_1),\dots]$. The properties of
$B$ are easy (compare, for example, \cite{UCA}, Ex.~8.5).

\begin{theorem} The principal ideal $m=tB\subset B$ is maximal, with
$B/m=k=A/(t)$, and the localisation $B_m$ is a DVR with the same parameter
$t$. \end{theorem}

\begin{pf} Consider the natural ``evaluation'' homomorphism $B\to k=A/(t)$
defined by $t\mapsto 0$, $z_r\mapsto a_r$. This is obviously surjective, and
by (\ref{eq:z_r=}), the kernel is the principal ideal $m=tB$. The localisation
$B_m$ is a local ring; its maximal ideal $mB_m=tB_m$ is principal; and
$\bigcap(t^n)=0$ in $B_m$, because $B_m\subset\hA$, and the same holds there.
This proves that $B$ is a DVR (see, for example, \cite{UCA}, Proposition~8.4)
with local parameter $t$, residue field $B/tB=A/(t)=k$ and $\Frac B=\Frac
A(z)$. \QED\end{pf}

Now the big one: set
 $$
C=A\bigl[t(z_0-a_0),\bigl\{(z_i-a_i)^2\bigr\}_{i=0}^\infty\bigr]\subset B.
 $$

\begin{theorem} $M=\bigl(t,t(z_0-a_0)\bigr)\subset C$ is a maximal ideal with
$C/M=k=A/(t)$, and the localisation $C_M$ has the following properties:
 \begin{enumerate}\renewcommand\labelenumi{(\roman{enumi})}
 \item $B$ and $C$ have the same field of fractions:
 $$
\Frac B=\Frac C=\Frac A(z).
 $$
 \item $B_m$ is integral over $C_M$, so that $B_m=\wave{C_M}$
(the normalisation).
 \item $C_M$ is a $1$-dimensional Noetherian local ring.
 \item $B_m$ is not finite as a $C_M$-module.
 \end{enumerate}

\end{theorem}

Statements (i) and (ii) are immediate. The surprise, of course, is that $C$ is
Noetherian.

\begin{pf} Manipulating the identities (\ref{eq:z_r=})--(\ref{eq:tnrz_r=})
gives two standard tricks. First, by (\ref{eq:tnrz_r=}), the difference
between $t(z_0-a_0)$ and $t^{n_r+1}(z_r-a_r)$ is an element of $A$ for any
$r\ge0$. This allows me to replace $t^{n_i+1}(z_i-a_i)$ wherever it appears by
an element of $A$ plus $t^{n_j+1}(z_j-a_j)$ with $j>i$.

Second, consider the identity
 \begin{equation}
\begin{aligned}
(z_{i-1}-a_{i-1})^2
&\hbox{}=(t^{m_i}z_i)^2\cr
&\hbox{}=t^{2m_i}((z_i-a_i)^2-a_i^2)+2a_it^{2m_i}z_i.
 \label{eq:(zi-ai)2}\end{aligned}
 \end{equation}
It's easy to check that both terms on the right are in $tC$: the second,
because of (\ref{eq:tnrz_r=}) and the assumption $2m_r\ge n_r+2$ (see
(\ref{ineq:nr}) and (\ref{m_r=})). A first consequence is that {\it the kernel
of the map $C\onto k$ defined by the evaluation $t\mapsto0$ and $z_i\mapsto
a_i$ is the maximal ideal $M=(t,t(z_0-a_0))$.}

The second trick allows me to replace $(z_{i-1}-a_{i-1})^2$ wherever it appears
by
 $$
t^{2m_i}(z_i-a_i)^2+\text{multiple of $t^{n_i+2}(z_i-a_i)$}+
\text{element of $A$}.
 $$
Performing these two tricks repeatedly gives that, {\em for any specified
$r\ge0$
and $N>0$, any element $f\in C$ can be written}
 \begin{equation}
f=X+Yt^{n_r+1}(z_r-a_r)+t^NZ,\quad\text{with $X,Y\in A$ and $Z\in C$.}
 \label{eq:f=X+Y}\end{equation}

 \begin{claim} For $0\ne f\in M$, the principal ideal $fC_M$ contains a power
of $t$. \end{claim}

\begin{pfof}{Claim} There exists $N$ such that $f\notin t^N\hA$. Choose $r$
with $n_r\ge N-1$, and consider the expression (\ref{eq:f=X+Y}). Then
necessarily
$X=t^nu$ with
$n<N$ and $u$ a unit of $A$. Dividing through by $u$, I assume that $X=t^n$,
and
 $$
f=t^n(1+t^{N-n}Z)+Yt^{n_r+1}(z_r-a_r).
 $$
To prove the claim, multiply $f$ by $g=t^n(1+t^{N-n}Z)-Yt^{n_r+1}(z_r-a_r)$:
 \begin{equation}
fg=t^{2n}(1+t^{N-n}Z)^2-Y^2t^{2n_r+2}(z_r-a_r)^2.
 \nonumber\end{equation}
This is obviously of the form $t^{2n}$ times an element of $C\setminus M$.
\QED \end{pfof}

I prove that the local ring $C_M$ is Noetherian and 1-dimensional. It is clear
from (\ref{eq:f=X+Y}) that $C_M/t^NC_M$ is generated over $A/(t^N)$ by 1 and
$t^{n_r+1}z_r$, and therefore is a Noetherian $A$-module. Now any nonzero
ideal $I\subset C_M$ contains $t^N$ for some $N$, and then the quotient ring
$C_M/I$ is also Noetherian. Therefore bigger ideals $I\subset J\subset C_M$
have the a.c.c. A nonzero prime ideal of $C_M$ contains some $t^N$, and
therefore also $t$ and $t(z_0-a_0)$, so that $\Spec C_M=\bigl\{0,MC_M\bigr\}$.

Under the assumption that $z$ is transcendental, I now prove that $B_m$ is not
finite over $C_M$, arguing by contradiction. Since $C_M$ is Noetherian, if
$B_m$ were finite, it would be a Noetherian $C_M$-module. Consider the
ascending chain of submodules generated by $\bigl\{(z_i-a_i)\bigr\}_{i<r}$;
for some $r$, I get a relation
 \begin{equation}
z_r-a_r=\sum_{i=0}^{r-1}g_i(z_i-a_i)\quad\text{with $g_i\in C_M$.}
 \label{eq:7}\end{equation}
Writing $g_i=f_i/f_r\in C_M$ gives
 \begin{equation}
f_r(z_r-a_r)=\sum_{i=0}^{r-1}f_i(z_i-a_i)
\quad\text{with $f_0,\dots,f_r\in C$ and $f_r\notin M$.}
 \label{eq:8}\end{equation}
Now multiplying (\ref{eq:8}) by $t^{n_r}$ and using (\ref{eq:tnrz_r=}) gives
 \begin{equation}
f_r(z-\sum_{j=0}^{r+1}a_jt^{n_j})=\sum_{i=0}^{r-1}f_it^{n_r-n_i}(z-\sum_{j=0}^{i+1}
a_jt^{n_j}).
 \label{eq:9}\end{equation}
Now all the $f_i\in C$ are polynomials in $z$ with coefficients in $A$, and
the left-hand side is a unit times $z$ (because $f_r\notin M$), whereas every
coefficient on the right-hand side is divisible by $t$. Therefore (\ref{eq:9})
is a nontrivial polynomial relation $F(z)=0$ with coefficients in $A$. This
contradiction completes the proof of the theorem. \QED \end{pf}

 \paragraph{Exercises}
 \begin{enumerate} \item Use (\ref{eq:(zi-ai)2}) to prove that $M^2=tM$.
 \item Prove that $t^{n_r}(z_r-a_r)\notin C$ for any $r\ge0$. [Hint: following
the method of (\ref{eq:9}), use $t^{n_r}(z_r-a_r)\in C$ to derive an algebraic
dependence relation for $z$ over $A$.]
 \end{enumerate}

 \paragraph{History} My treatment follows Akizuki in all essentials. Clearly
under the influence of the papers of Krull and his followers, Akizuki only
considers the case where $\hA=\Z_p$ is the ring of $p$-adic integers. His
proof that $C$ is infinite over $B$ is indirect. He argues by contradiction,
based on the notion of ``analytically unramified'' (in later terminology): the
element $x_r=t^{n_r+1}(z_r-a_r)\notin tC$ (by Ex.~2 above), but $x_r^2\in
t^{2n_r+2}C$. Thus $x=\lim_{r\to\infty}x_r$ is a nilpotent element of the
$t$-adic completion of $C$.

As discussed in \cite{UCA}, 9.4, the real point of this counterexample, and of
those of Nagata (see the appendix to \cite{Nagata}) is that there is really no
hope of making everything that works for geometric rings go through for
Noetherian rings. At some time you have to make assumptions of a concrete
nature, for example that your ring is finitely generated over $k$ or $\Z$.

 \paragraph{Geometric interpretation of $B$} If $A=\C[t]_{(0)}$ and the power
series $z$ has positive radius of convergence, I can consider the analytic arc
$\Ga\subset\C^2$ defined by $(z=z(t))$. There is an obvious sense in which
$B_m$ is the ring of regular functions on $\Ga$ that are restrictions of
rational functions of $t,z$.

More algebraically, for each $r$, I can view $\aff_r=\Spec A[z_r]$ as the
``affine plane'' with coordinates $t,z_r$, or its germ at $(t=0,z=a_r)$. The
inclusion of rings $\Spec A[z_{r-1}]\subset\Spec A[z_r]$ corresponds to the
``blow-up'' $\aff_r\to\aff_{r-1}$ defined by $(t,z_r)\mapsto
(t,a_r+t^{n_r}z_{r-1})$. The limit $\Spec B$ is the surface in infinite
dimensional space defined by the relations (\ref{eq:tnrz_r=}). The projection
to each $\aff_r$ can be viewed as an infinitely thin cusp-shaped region around
the analytic arc $z=z(t)$.

Rings like $B$ are interesting because of their proclivity to dimensional
ambi\-guity, arising from the question as to whether or not (\ref{z=z(t)}) is
a functional dependence relation $z=z(t)$. This ambiguity is the starting point
for Nagata's examples of noncatenary rings, see \cite{Nagata}, Example~2,
p.~203 or \cite{UCA}, \hbox{Example~9.4}, (2). As we have seen, $B$ becomes a
DVR when localised at $(t)$, because modulo $t^N$ the identities
(\ref{eq:z_r=})--(\ref{eq:tnrz_r=}) imply the ``obvious'' relation $z=z(t)$.
On the other hand, if I delocalise $A$ (taking, say, $A=k[t]$), I can pass to
the ring of fractions $B[1/t]$. Then the identities (\ref{eq:z_r=}), give all
the $z_i$ as functions of $z=z_0$, so that, assuming $z$ is transcendental over
$A$, $B[1/t]=A[1/t][z]$ is clearly 2-dimensional (for example, $B=k[t,z][1/t]$
is just polynomial functions on the $z,t$ plane).

 \paragraph{Geometric interpretation of $C$} Even after the event, I don't
know how to motivate Akizuki's example to make it completely natural, and it's
hard to imagine how he discovered his incredibly ingenious construction.

For what it's worth, I have in mind the following geometric picture, by
analogy with the above picture of $B$: the ring $C\subset B$ is the union over
$r$ of subrings
 $$
A[t^{n_r+1}z_r,(z_r-a_r)^2]\subset A[z_r].
 $$
In other words, the monomials that are missing are $t^iz_r$ for
$i=1,\dots,n_r$. This can be interpreted as creasing the $z,t$ plane $\aff_r$
along the analytic arc $\Ga$ to have a cusp for $n_r+1$ infinitesimal steps.
Of course, it's hard to predict on the basis of the geometric picture why such
a weird procedure should lead to a Noetherian ring.

 \paragraph{Thanks} To John Moody and Shigeru Mukai for helpful discussions.

 \paragraph{Ad} This note is a free sample of my forthcoming book \cite{UCA}.
Place your order soon to avoid disappointment.


\begin{thebibliography}{Akizuki}

 \bibitem[Akizuki]{Akizuki} Y. Akizuki, Einige Bemerkungen \"uber prim\"are
Integrit\"atsbereiche mit Teilerkettensatz, Proc.\ Phys.--Math.\ Soc.\ Japan,
{\bf 17} (1935), pp.~327--336

 \bibitem[Nagata]{Nagata} M. Nagata, Local rings, John Wiley Interscience, New
York, London, 1962

 \bibitem[UCA]{UCA} M. Reid, Undergraduate commutative algebra,
C.U.P., Cambridge 1995

 \end{thebibliography}
\end{document}